# SCAD-PENALIZED REGRESSION IN HIGH-DIMENSIONAL PARTIALLY LINEAR MODELS

By Huiliang Xie and Jian Huang

*University of Miami and University of Iowa*

We consider the problem of simultaneous variable selection and estimation in partially linear models with a divergent number of covariates in the linear part, under the assumption that the vector of regression coefficients is sparse. We apply the SCAD penalty to achieve sparsity in the linear part and use polynomial splines to estimate the nonparametric component. Under reasonable conditions, it is shown that consistency in terms of variable selection and estimation can be achieved simultaneously for the linear and nonparametric components. Furthermore, the SCAD-penalized estimators of the nonzero coefficients are shown to have the asymptotic oracle property, in the sense that it is asymptotically normal with the same means and covariances that they would have if the zero coefficients were known in advance. The finite sample behavior of the SCAD-penalized estimators is evaluated with simulation and illustrated with a data set.

**1. Introduction.** Consider a partially linear model (PLM)

$$Y = \mathbf{X}'\boldsymbol{\beta} + g(T) + \varepsilon,$$

where $\boldsymbol{\beta}$ is a $p \times 1$ vector of regression coefficients associated with $\mathbf{X}$, and $g$ is an unknown function of $T$. In this model, the mean response is linearly related to $\mathbf{X}$, while its relation with $T$ is not specified up to any finite number of parameters. This model combines the flexibility of nonparametric regression and parsimony of linear regression. When the relation between $Y$ and $\mathbf{X}$ is of main interest and can be approximated by a linear function, it offers more interpretability than a purely nonparametric model.

We consider the problem of simultaneous variable selection and estimation in the PLM when $p$ is large, in the sense that $p \to \infty$ as the sample size $n \to \infty$. For finite-dimensional $\boldsymbol{\beta}$, several approaches have been proposed to estimate $\boldsymbol{\beta}$ and $g$. Examples include the partial spline estimator [Wahba (1984), Engle et al. (1986) and Heckman (1986)] and the partial residual estimator [Robinson (1988), Speckman (1988) and Chen (1988)]. Under appropriate assumptions about the smoothness of $g$ and the structure of $\mathbf{X}$, these estimators of $\boldsymbol{\beta}$ were shown to be $\sqrt{n}$-consistent and asymptotically normal. It was also shown that the estimators of $g$ can converge at the optimal rate in the purely nonparametric regression determined in Stone (1980, 1982). Fan and Li (2004) considered variable selection in the semiparametric models in the context of longitudinal data analysis, assuming a framework with a fixed set of covariates as $n$ increases. In these studies, either the dimension of the covariate vector $\mathbf{X}$ was fixed or the problem of variable selection in $\mathbf{X}$ via penalization was not considered. However, the results for the PLM with a finite-dimensional $\boldsymbol{\beta}$ and those for the semiparametric models in general are not applicable to the PLM with a divergent number of covariates. Indeed, it appears that there is no systematic theoretical investigation of estimation in semiparametric models with a high-dimensional parametric component.

We are particularly interested in $\boldsymbol{\beta}$ when it is sparse, in the sense that many of its elements are zero. Our work is motivated by biomedical studies that investigate the relationship between a phenotype of interest and genomic measurements such as microarray data. In many such studies, in addition to genomic measurements, other types of measurements, such as clinical or environmental covariates, are also available. To obtain unbiased estimates of genomic effects, it is necessary to take into account these covariates. Assuming a sparse model is often reasonable with genomic data. This is because, although the total number of measurements can be large, the number of important ones is usually relatively small. In these problems, selection of important covariates is often one of the most important goals in the analysis. The $p \to \infty$ framework allows us to address the concerns as to how the nonparametric term is going to affect the estimation and variable selection of $\boldsymbol{\beta}$, and whether the rate at which the nonparametric estimator converges can be maintained with a divergent $p$.

We use the SCAD method to achieve simultaneous consistent variable selection and estimation of $\boldsymbol{\beta}$. The SCAD method is proposed by Fan and Li (2001) in a general parametric framework for variable selection and efficient estimation. This method uses a specially designed penalty function, the smoothly clipped absolute deviation (hence the name SCAD), as adopted in Fan and Li (2004). We estimate the nonparametric component $g$ using the partial residual method with the B-spline bases. The resulting estimator of $\boldsymbol{\beta}$ maintains the oracle property of the SCAD-penalized estimators in parametric settings. Here, the oracle property means that the estimator can correctly select the nonzero coefficients with probability converging to one, and that the estimators of the nonzero coefficients are asymptotically normal with the same means and covariances that they would have if the zero coefficients were known in advance. Therefore, an oracle estimator is asymptotically as efficient as the ideal estimator assisted by an oracle who knows which coefficients are nonzero. Meanwhile, convergence of the estimator of $g$ in the SCAD-penalized partially linear regression still reaches the optimal global rate.

Investigations on the asymptotic properties of penalized estimation in parametric models when the number of covariates is fixed include Knight and Fu (2000) and Fan and Li (2001). Fan and Peng (2004) considered the same problem when the number of parameters diverges, where they showed that there exist local maximizers of the penalized likelihood that have an oracle property. Huang, Horowitz and Ma (2008) studied the bridge estimators with a divergent number of covariates in a linear regression model and showed that the bridge estimators have an oracle property if the bridge index is strictly between 0 and 1. Several recent studies have considered the asymptotic properties of the LASSO method in high-dimensional settings. Examples include: Meinshausen and Buhlmann (2006), van de Geer (2008), Zhang and Huang (2008) and Zhao and Yu (2006). In these studies, the convexity property of the LASSO penalty is critical to the results. However, since the SCAD penalty is not convex, the methods that utilize convexity are not applicable in the present setting. Furthermore, the PLM models we consider here are semiparametric. The asymptotic analysis of such semiparametric models in high-dimensional settings appears to be considerably more complicated than those in the linear regression models.

The rest of this article is organized as follows. In Section 2, we define the SCAD-penalized estimator $(\widehat{\boldsymbol{\beta}}_n, \widehat{g}_n)$ in the PLM, abbreviated as SCAD-PLM estimator hereafter. The main results for the SCAD-PLM estimator are given in Section 3, including the consistency and oracle property of $\widehat{\boldsymbol{\beta}}_n$, as well as the rate of convergence of $\widehat{g}_n$. Section 4 deals with computing the PLM-SCAD estimator. The finite sample behavior of this estimator is illustrated with simulation studies and a real data example in Section 5. Extensions and concluding remarks are given in Section 6. The proofs are relegated to the Appendix.

**2. Penalized estimation in PLM with the SCAD penalty.** To make it explicit that the covariates and regression coefficients depend on $n$, we write the PLM

$$Y_i = \mathbf{X}_i^{(n)\prime} \boldsymbol{\beta}^{(n)} + g(T_i) + \varepsilon_i, \qquad i = 1, \ldots, n,$$

where $(\mathbf{X}_i^{(n)}, T_i, Y_i)$ are independent and identically distributed as $(\mathbf{X}^{(n)}, T, Y)$, and $\varepsilon_i$ is independent of $(\mathbf{X}_i^{(n)}, T_i)$, with mean 0 and variance $\sigma^2$. We assume that $T$ takes values in a compact interval, and, for simplicity, we assume this interval to be $[0,1]$. Let $\mathbf{Y} = (Y_1, \ldots, Y_n)'$, and let $\mathbb{X}^{(n)} = (X_{ij}, 1 \leq i \leq n, 1 \leq j \leq p_n)$ be the $n \times p_n$ design matrix associated with $\boldsymbol{\beta}^{(n)}$. In sparse models, the $p_n$ covariates can be classified into two categories: the important ones whose corresponding coefficients are nonzero and the trivial ones that actually are not present in the underlying model. For convenience of notation, we write

$$\boldsymbol{\beta}^{(n)} = (\boldsymbol{\beta}_1^{(n)\prime}, \boldsymbol{\beta}_2^{(n)\prime})', \tag{1}$$

where $\boldsymbol{\beta}_1^{(n)\prime} = (\beta_1^{(n)}, \ldots, \beta_{k_n}^{(n)})$ and $\boldsymbol{\beta}_2^{(n)\prime} = (0, \ldots, 0)$. Here $k_n (\leq p_n)$ is the number of nontrivial covariates. Let $m_n = p_n - k_n$ be the number of zero coefficients.

We use the polynomial splines to approximate $g$. For a positive integer $M_n$, let $\Delta_n = \{\xi_{n\nu}\}_{\nu=1}^{M_n}$ be a partition of $[0,1]$ into $M_n + 1$ subintervals $I_{n\nu} = [\xi_{n\nu}, \xi_{n,\nu+1}) : \nu = 0, \ldots, M_n - 1$ and $I_{nM_n} = [\xi_{nM_n}, 1]$. Here, $\xi_{n0} = 0$ and $\xi_{n,M_n+1} = 1$. Denote the largest mesh size of $\Delta_n$, $\max_{0 \leq \nu \leq M_n}\{\xi_{n,\nu+1} - \xi_{n\nu}\}$, by $\overline{\Delta}_n$. Throughout the article, we assume $\overline{\Delta}_n = O(M_n^{-1})$. Let $\mathcal{S}_m(\Delta_n)$ be the space of polynomial splines of order $m$ with simple knots at the points $\xi_{n1}, \ldots, \xi_{nM_n}$. This space consists of all functions $s$ with these two properties:

(i) restricted to any interval $I_{n\nu}$ $(0 \leq \nu \leq M_n)$, $s$ is a polynomial of order $m$;
(ii) if $m \geq 2$, $s$ is $m-2$ times continuously differentiable on $[0,1]$.

According to Corollary 4.10 in Schumaker (1981), there is a local basis $\{B_{nw}, 1 \leq w \leq q_n\}$ for $\mathcal{S}_m(\Delta_n)$, where $q_n = M_n + m$ is the dimension of $\mathcal{S}_m(\Delta_n)$. Let

$$\mathbf{Z}(t; \Delta_n)' = (B_{n1}(t), \ldots, B_{nq_n}(t))$$

and $\mathbb{Z}^{(n)}$ be the $n \times q_n$ matrix whose $i$th row is $\mathbf{Z}(T_i; \Delta_n)'$. Any $s \in \mathcal{S}_m(\Delta_n)$ can be written $s(t) = \mathbf{Z}(t; \Delta_n)' \mathbf{a}^{(n)}$ for a $q_n \times 1$ vector $\mathbf{a}^{(n)}$. We try to find the $s$ in $\mathcal{S}_m(\Delta_n)$ that is close to $g$. Under reasonable smoothness conditions, $g$ can be well approximated by elements in $\mathcal{S}$. Thus, the problem of estimating $g$ becomes that of estimating $\mathbf{a}^{(n)}$.

Given $a > 2$ and $\lambda > 0$, the SCAD penalty at $\theta$ is

$$p_\lambda(\theta; a) = \begin{cases} \lambda|\theta|, & |\theta| \leq \lambda, \\ -(\theta^2 - 2a\lambda|\theta| + \lambda^2)/[2(a-1)], & \lambda < |\theta| \leq a\lambda, \\ (a+1)\lambda^2/2, & |\theta| > a\lambda. \end{cases}$$

The SCAD penalty is continuously differentiable on $(-\infty, 0) \cup (0, \infty)$ but singular at 0. Its derivative vanishes outside $[-a\lambda, a\lambda]$. As a consequence, SCAD penalized regression can produce sparse solutions and unbiased estimates for large coefficients. More details of the penalty can be found in Fan and Li (2001).

The penalized least squares objective function for estimating $\boldsymbol{\beta}^{(n)}$ and $\mathbf{a}^{(n)}$ with the SCAD penalty is

$$
\begin{aligned}
Q_n(\mathbf{b}^{(n)}, \mathbf{a}^{(n)}; \lambda_n, a, \Delta_n, m) & \\
(2) \qquad & = \|\mathbf{Y} - \mathbb{X}^{(n)}\mathbf{b}^{(n)} - \mathbb{Z}^{(n)}\mathbf{a}^{(n)}\|^2 + n \sum_{j=1}^{p_n} p_{\lambda_n}(b_j^{(n)}; a).
\end{aligned}
$$

Let

$$(\widehat{\boldsymbol{\beta}}_n^{(n)}, \widehat{\boldsymbol{\alpha}}_n^{(n)}) = \arg\min Q_n(\mathbf{b}^{(n)}, \mathbf{a}^{(n)}; \lambda_n, a, \Delta_n, m).$$

The SCAD-PLM estimators of $\boldsymbol{\beta}$ and $g$ are $\widehat{\boldsymbol{\beta}}_n$ and $\widehat{g}_n(t) \equiv \mathbf{Z}(t; \Delta_n)' \widehat{\boldsymbol{\alpha}}_n^{(n)}$, respectively.

The polynomial splines were also used by Huang (1999) in the partially linear Cox models. Some computational conveniences were also discussed there. We limit our search for the estimate of $g$ to the space of polynomial splines of order $m$, instead of the larger space of piecewise polynomials of order $m$, with the goal to find a smooth estimator of $g$. Unlike the basis pursuit in nonparametric regression, no penalty is imposed on the estimator of the nonparametric part, as our interest lies in the variable selection with regard to the parametric part.

For any $\mathbf{b}^{(n)}$, the $\mathbf{a}^{(n)}$ that minimizes $Q_n$ necessarily satisfies

$$\mathbb{Z}^{(n)\prime}\mathbb{Z}^{(n)}\mathbf{a}^{(n)} = \mathbb{Z}^{(n)\prime}(\mathbf{Y} - \mathbb{X}^{(n)\prime}\mathbf{b}^{(n)}).$$

Let $P_{\mathbf{Z}}^{(n)} = \mathbb{Z}^{(n)}(\mathbb{Z}^{(n)\prime}\mathbb{Z}^{(n)})^{-1}\mathbb{Z}^{(n)\prime}$ be the projection matrix of the column space of $\mathbb{Z}^{(n)}$. The profile objective function of the parametric part becomes

$$
\widetilde{Q}_n(\mathbf{b}^{(n)}; \lambda_n, a, \Delta_n, m) = \|(I - P_{\mathbf{Z}}^{(n)})(\mathbf{Y} - \mathbb{X}^{(n)}\mathbf{b}^{(n)})\|^2 + n \sum_{j=1}^{p_n} p_{\lambda_n}(b_j^{(n)}; a).
$$

(3)

Then, $\widehat{\boldsymbol{\beta}}_n^{(n)} = \arg\min \widetilde{Q}_n(\mathbf{b}^{(n)}; \lambda_n, a, \Delta_n, m)$. Because the profile objective function does not involve $\mathbf{a}^{(n)}$ and has an explicit form, it is useful for both theoretical investigation and computation. We will use it to established the asymptotic properties of $\widehat{\boldsymbol{\beta}}_n^{(n)}$. Computationally, this expression can be used to first obtain $\widehat{\boldsymbol{\beta}}_n^{(n)}$. Then, $\widehat{\boldsymbol{\alpha}}_n^{(n)}$ can be computed using the resulting residuals as the response for the covariate matrix $\mathbb{Z}^{(n)}$.

**3. Asymptotic properties of the PLM-SCAD estimator.** In this section we state the results of the asymptotic properties of the PLM-SCAD estimator. First, we define some notation. Let $\theta_j^{(n)}(t) = E[X_j^{(n)}|T=t]$ for $j = 1, \ldots, p_n$. Let the $p_n \times p_n$ conditional variance-covariance matrix of $(\mathbf{X}^{(n)}|T=$

$t$) be $\Sigma^{(n)}(t)$. Let $\mathbf{e}^{(n)} = \mathbf{X}^{(n)} - E[\mathbf{X}^{(n)}|T]$. We can write $\Sigma^{(n)}(t) = \mathrm{Var}(\mathbf{e}^{(n)}|T = t)$. Denote the unconditional variance-covariance matrix of $\mathbf{e}^{(n)}$ by $\Xi^{(n)}$. We have $\Xi^{(n)} = E[\Sigma^{(n)}(T)]$. We assume the following conditions on the smoothness of $g$ and $\theta_j^{(n)}, 1 \leq j \leq p_n$.

CONDITION 1. There are absolute constants $\gamma_\theta > 0$ and $M_\theta > 0$, such that
$$\sup_{n \geq 1} \sup_{1 \leq j \leq p_n} |\theta_{nj}^{(r_\theta)}(t_2) - \theta_{nj}^{(r_\theta)}(t_1)| \leq M_\theta |t_2 - t_1|^{\gamma_\theta} \quad \text{for } 0 \leq t_1, t_2 \leq 1,$$
and the degree of the polynomial spline $m - 1 \geq r_\theta$. Let $s_\theta = r_\theta + \gamma_\theta$.

CONDITION 2. There exists an absolute constant $\sigma_{4e}$, such that, for all $n$ and $1 \leq j \leq p_n$,
$$E[e_j^{(n)4}|T] \leq \sigma_{4e} \quad \text{almost surely.}$$

CONDITION 3. There are absolute constants $\gamma_g > 0$ and $M_g > 0$, such that
$$|g^{(r_g)}(t_2) - g^{(r_g)}(t_1)| \leq M_g |t_2 - t_1|^{\gamma_g} \quad \text{for } 0 \leq t_1, t_2 \leq 1,$$
with $r_g \leq m - 1$. Let $s_g = r_g + \gamma_g$.

As in nonparametric regression, we allow $M_n \to \infty$, but $M_n = o(n)$. In addition, we assume that the tuning parameter $\lambda_n \to 0$ as $n \to \infty$. This is the assumption adopted in nonconcave penalized regression [Fan and Peng (2004)]. For convenience, all the other conditions required for the conclusions in this section are listed here.

(A1) (a) $\lim_{n \to \infty} p_n^2/n = 0$; (b) $\lim_{n \to \infty} p_n^2 M_n^2/n^2 = 0$; (c) $\lim_{n \to \infty} p_n/M_n^{s_\theta} = 0$.

(A2) The smallest eigenvalue of $\Xi^{(n)}$, denoted by $\lambda_{\min}(\Xi^{(n)})$, satisfies
$$\liminf_{n \to \infty} \lambda_{\min}(\Xi^{(n)}) = c_\lambda > 0.$$

(A3) $\lambda_n = o(k_n^{-1/2})$.

(A4) $\liminf_{n \to \infty} \min_{1 \leq j \leq k_n} |\beta_j^{(n)}| = c_\beta > 0$.

(A5) Let $\lambda_{\max}(\Xi^{(n)})$ be the largest eigenvalue of $\Xi^{(n)}$. (a) $\lim \sqrt{p_n \lambda_{\max}(\Xi^{(n)})}/(\sqrt{n}\lambda_n) = 0$; (b) $\lim \sqrt{p_n \lambda_{\max}(\Xi^{(n)})}/(M_n^{s_g} \lambda_n) = 0$.

(A6) Suppose for all $t$ in $[0,1]$, $\mathrm{tr}(\Sigma_{11}^{(n)}(t)) \leq \mathrm{tr}(\Sigma_{u,11}^{(n)})$ and the latter satisfies $\lim \sqrt{\mathrm{tr}(\Sigma_{u,11}^{(n)})} M_n^{-s_g} = 0$ and $\lim \mathrm{tr}(\Sigma_{u,11}^{(n)}) M_n/n = 0$.

(A7) $\lim \sqrt{n} M_n^{-(s_g + s_\theta)} = 0$.

THEOREM 1 (Consistency of $\widehat{\boldsymbol{\beta}}^{(n)}$). *Under* (A1)–(A2),
$$\|\widehat{\boldsymbol{\beta}}^{(n)} - \boldsymbol{\beta}^{(n)}\| = O_P(\sqrt{p_n/n} + M_n^{-s_g} + \sqrt{k_n}\lambda_n).$$

*Thus, under* (A1)–(A3), $\|\widehat{\boldsymbol{\beta}}^{(n)} - \boldsymbol{\beta}^{(n)}\| \xrightarrow{P} 0.$

This theorem establishes the consistency of the PLM-SCAD estimator of the parametric part, without the local restriction in Theorem 1 of Fan and Peng (2004) and Theorem 2 of Fan and Li (2004). (A1) requires the number of covariates considered not to increase at rates faster than $\sqrt{n}$ and $M_n^{1/s_\theta}$. (A2) is a requirement for model identifiability. It assumes that $\Xi^{(n)}$ is positive definite, so that no random variable of the form $\sum_{j=1}^{p_n} c_j X_j^{(n)}$, where $c_j$'s are constants, can be functionally related to $T$. When $p_n$ increases with $n$, $\Xi^{(n)}$ needs to be bounded away from any singular matrix. The assumption about $\lambda_n$, (A3), says that $\lambda_n$ should converge to 0 fast enough so that the penalty would not introduce any bias. The rate at which $\lambda_n$ goes to 0 only depends on $k_n$. It is interesting to note that the smoothness index $s_g$ of $g$ and the number of spline bases $M_n$ affects the rate of convergence of $\widehat{\boldsymbol{\beta}}^{(n)}$ by contributing a term $M_n^{-s_g}$. When $p_n$ is bounded and no SCAD penalty is imposed ($\lambda_n = 0$), the convergence rate is $O(n^{-1/2} + M_n^{-s_g})$, which is consistent with Theorem 2 of Chen (1988).

Corresponding to the partition in (1), write $\widehat{\boldsymbol{\beta}}^{(n)} = (\widehat{\boldsymbol{\beta}}_1^{(n)\prime}, \widehat{\boldsymbol{\beta}}_2^{(n)\prime})'$, where $\widehat{\boldsymbol{\beta}}_1^{(n)\prime}$ and $\widehat{\boldsymbol{\beta}}_2^{(n)\prime}$ are vectors of length $k_n$ and $m_n$, respectively. The theorem below shows that all the covariates with zero coefficients can be detected simultaneously with probability tending to 1, provided that $\lambda_n$ does not converge to 0 too fast.

THEOREM 2 (Variable selection in $\mathbf{X}^{(n)}$). *Assume all the $e_j^{(n)}$'s support sets are contained in a compact set in $\mathcal{R}$. Under* (A1)–(A5), $\lim_{n\to\infty} P(\widehat{\boldsymbol{\beta}}_2^{(n)} = \mathbf{0}) = 1$.

(A5) puts restriction on the largest eigenvalue of $\Xi^{(n)}$. In general, $\lambda_{\max}(\Xi^{(n)}) = O(p_n)$, as can be seen from
$$\lambda_{\max}(\Xi^{(n)}) < \operatorname{tr}(\Xi^{(n)}) \leq p_n\sqrt{\sigma_{4e}}.$$

There is the question of whether there exists a $\lambda_n$ that satisfies both (A3) and (A5). It can be checked that, if $p_n = o(n^{1/3})$, there exists $\lambda_n$, such that (A3) and (A5) hold. When $k_n$ is bounded, the existence of such $\lambda_n$ only requires that $p_n = o(n^{1/2})$. This relaxation also holds for the case when $\lambda_{\max}(\Xi^{(n)})$ is bounded from above.

By Theorem 2, $\widehat{\boldsymbol{\beta}}_2^{(n)}$ degenerates at $\mathbf{0}_{m_n}$, with probability converging to 1. We now consider the asymptotic distribution of $\widehat{\boldsymbol{\beta}}_1^{(n)}$. According to the partition of $\boldsymbol{\beta}^{(n)}$ in (1), write $\mathbb{X}^{(n)}$ and $\Xi^{(n)}$ in the block form:

$$\mathbb{X}^{(n)} = (\underbrace{\mathbb{X}_1^{(n)}}_{n \times k_n} \underbrace{\mathbb{X}_2^{(n)}}_{n \times m_n}), \qquad \Xi^{(n)} = \begin{matrix} & k_n & m_n \\ k_n & \\ m_n & \end{matrix} \begin{pmatrix} \Xi_{11}^{(n)} & \Xi_{12}^{(n)} \\ \Xi_{21}^{(n)} & \Xi_{22}^{(n)} \end{pmatrix}.$$

Let $\mathbf{A}_n$ be a nonrandom $\iota \times k_n$ matrix with full row rank, and

$$\Sigma_n = n^2 \mathbf{A}_n [\mathbb{X}_1^{(n)\prime}(I - P_{\mathbf{Z}}^{(n)})\mathbb{X}_1^{(n)}]^{-1} \Xi_{11}^{(n)} [\mathbb{X}_1^{(n)\prime}(I - P_{\mathbf{Z}}^{(n)})\mathbb{X}_1^{(n)}]^{-1} \mathbf{A}_n'.$$

THEOREM 3 (Asymptotic distribution of $\widehat{\boldsymbol{\beta}}^{(n)}$). *Suppose that all the support sets of $e_j^{(n)}$'s are contained in a compact set in $\mathcal{R}, j = 1, \ldots, p_n$. Then, under* (A1)–(A7),

(4) $$\sqrt{n} \Sigma_n^{-1/2} \mathbf{A}_n (\widehat{\boldsymbol{\beta}}_1^{(n)} - \boldsymbol{\beta}_1^{(n)}) \xrightarrow{d} N(\mathbf{0}_\iota, \sigma^2 \mathbf{I}_\iota).$$

The asymptotic distribution result can be used to construct asymptotic confidence intervals for any fixed number of coefficients simultaneously.

In (4), we used the inverse of $\mathbb{X}_1^{(n)\prime}(I - P_{\mathbf{Z}}^{(n)})\mathbb{X}_1^{(n)}$ and that of $\Sigma_n$. Under assumption (A2), by Theorem 4.3.1 in Wang and Jia (1993), the smallest eigenvalue of $\Xi_{11}^{(n)}$ is no less than $c_\lambda$ and bounded away from 0. By Lemma 1 in the Appendix, $\mathbb{X}_1^{(n)\prime}(I - P_{\mathbf{Z}}^{(n)})\mathbb{X}_1^{(n)}$ is invertible with probability tending to 1. The invertibility of $\Sigma_n$ then follows from the full row rank restriction on $\mathbf{A}_n$.

(A6) may appear a little abrupt. It requires $\sum_{j=1}^{k_n} \text{Var}(e_j^{(n)}|T=t)$ to be less than the trace of a $k_n \times k_n$ matrix $\Sigma_{u,11}^{(n)}$ as $t$ ranges over $[0,1]$, which is considerably weaker than the assumption that $\Sigma_{u,11}^{(n)} - \Sigma_{11}^{(n)}(t)$ is a nonnegative definite matrix for any $t \in [0,1]$. We can also replace $\text{tr}(\Sigma_{u,11}^{(n)})$ by $k_n$ in the assumption, since for all $t$, $\sum_{j=1}^{k_n} \text{Var}(e_j^{(n)}|T=t) \leq k_n \sqrt{C_e}$. (A7) requires that $g$ and $\theta_j^{(n)}$ be smooth enough. Intuitively, a smooth $g$ makes it easier to estimate $\boldsymbol{\beta}$. The smoothness requirement on $\theta_j^{(n)}$ also makes sense, since this helps to remove effect of $T$ on $X_j^{(n)}$, and the estimation of $\boldsymbol{\beta}$ is based on the relationship

$$Y - E[Y|T] = (\mathbf{X} - E[\mathbf{X}|T])\beta + \varepsilon.$$

We now consider the consistency of $\widehat{g}_n$. Suppose that $T$ is an absolutely continuous random variable on $[0,1]$ with density $f_T$. We use the $L_2$ distance

$$\|\widehat{g}_n - g\|_T = \left\{ \int_0^1 [\widehat{g}_n(t) - g(t)]^2 f_T(t)\, dt \right\}^{1/2}.$$

This is the measure of distance between two functions that were used in Stone (1982, 1985). If our interest is confined to the estimation of $\boldsymbol{\beta}^{(n)}$, we should choose large $M_n$, unless computing comes into consideration. However, too large an $M_n$ would introduce too much variation and is detrimental to the estimation of $g$.

THEOREM 4 (Rate of convergence of $\widehat{g}_n$). *Suppose $M_n = o(\sqrt{n})$, $f_T(t)$ is bounded away from $0$ and infinity on $[0,1]$ and $E[\varepsilon^4] < \infty$. Under* (A1)–(A5),

$$\|\widehat{g}_n - g\|_T = O_P(k_n/\sqrt{n} + \sqrt{M_n/n} + \sqrt{k_n} M_n^{-s_g}).$$

In the special case of bounded $k_n$, Theorem 4 simplifies to the well-known result in nonparametric regression: $\|\widehat{g}_n - g\|_T = O_P(\sqrt{M_n/n} + M_n^{-s_g})$. When $M_n \sim n^{-1/(2s_g+1)}$, the convergence rate is optimal. However, the feasibility of such a choice requires $s_g > 1/2$. To have the asymptotic normality of $\widehat{\boldsymbol{\beta}}_1^{(n)}$ hold simultaneously, we also need $s_\theta > 1/2$. In the diverging $k_n$ case, the rate of convergence is determined by $k_n$, $p_n$, $M_n$, $s_g$ and $s_\theta$ jointly. With appropriate $s_g$, $s_\theta$ and $p_n$, the rate of convergence can be $n^{-1/2} k_n + k_n^{1/(4s_g+2)} n^{-s_g/(2s_g+1)}$.

**4. Computation.** The computation of the PLM-SCAD estimator involves the choice of $\lambda_n$. We first consider the estimation, as well as the standard error approximation of the estimator with a given $\lambda_n$, and then describe the generalized cross validation approach to choose appropriate $\lambda_n$ in the PLM.

The computation of $(\widehat{\boldsymbol{\beta}}^{(n)}, \widehat{g}_n)$ requires the minimization of (2). The projection approach adopted here converts this problem to the minimization of (3). In particular, given $m$ and a partition $\Delta_n$, a basis of $\mathcal{S}_m(\Delta_n)$ is given by $(B_{n1}, \ldots, B_{nq_n})$. The basis functions are evaluated at $T_i, i = 1, \ldots, n$, and form $\mathbb{Z}_n$. In Splus or R, this can be realized with the bs function. Regress each column of $\mathbb{X}^{(n)}$ and $\mathbf{Y}$ on $\mathbb{Z}_n$ separately. Denote the residuals by $\widetilde{\mathbb{X}}^{(n)}$ and $\widetilde{\mathbf{Y}}$. The minimization of (3) is now a nonconcave penalized regression problem, with observations $(\widetilde{\mathbb{X}}^{(n)}, \widetilde{\mathbf{Y}})$. So, the minorize–maximize (MM) algorithm described in Hunter and Li (2005) can be used to compute $\widehat{\boldsymbol{\beta}}^{(n)}$. We also standardize the columns of $\widetilde{\mathbb{X}}^{(n)}$, so the covariates with smaller variations will not be discriminated against. Once we have computed $\widehat{\boldsymbol{\beta}}^{(n)}$, the value of $g$ at any $t \in [0,1]$ is estimated by $\widehat{g}_n(t) = \mathbf{Z}(t;\Delta_n)'(\mathbb{Z}^{(n)\prime}\mathbb{Z}^{(n)})^{-1}\mathbb{Z}^{(n)\prime}(\mathbf{Y} - \mathbb{X}^{(n)}\widehat{\boldsymbol{\beta}}^{(n)})$.

The standard errors of the nonzero components of $\widehat{\boldsymbol{\beta}}^{(n)}$ can be derived from the Hessian matrix. For details, see Hunter and Li (2005) or Fan and Li (2001).

We choose $\lambda_n$ by minimizing the generalized cross validation score [Wahba (1990)] and fix $a = 3.7$, as suggested by Fan and Li (2001). Our preference of GCV over CV stems from as much its computation advantage as its comparable performance to CV in model selection, which have been discussed in Tibshirani (1996) and Kim, Kim and Kim (2006).

Note that here we use fixed partition $\Delta_n$ and $m$ in estimating the nonparametric component $g$. Data-driven choice of them may be desirable, which inevitably requires a good estimator of $\boldsymbol{\beta}^{(n)}$. In our simulations, $m = 4$ (cubic splines) and $M_n \leq 3$ with even partition of $[0,1]$ serves the purpose well.

**5. Numerical studies.** In this section, we illustrate the PLM-SCAD estimator's finite sample properties with examples. Examples 1 is a simulated example, and Example 2 explores a real data set. Throughout, we use $m = 4$, $M_n = 3$ and the sample quantiles of $T_i$'s as the knots.

EXAMPLE 1. In this study, we simulate $n = 100$ points $T_i, i = 1, \ldots, 100$, from the uniform distribution on $[0,1]$. For each $i$, $e_{ij}$'s are simulated to be normally distributed with autocorrelated variance structure $AR(\rho)$, such that

$$\text{Cov}(e_{ij}, e_{il}) = \rho^{|j-l|}, \qquad 1 \leq j, l \leq 10.$$

$X_{ij}$'s are then formed as follows:

$$X_{i1} = \sin(2T_i) + e_{i1}, \qquad X_{i2} = (0.5 + T_i)^{-2} + e_{i2},$$
$$X_{i3} = \exp(T_i) + e_{i3}, \qquad X_{i5} = (T_i - 0.7)^4 + e_{i5},$$
$$X_{i6} = T_i(1 + T_i^2)^{-1} + e_{i6}, \qquad X_{i7} = \sqrt{1 + T_i} + e_{i7},$$
$$X_{i8} = \log(3T_i + 8) + e_{i8}, \qquad X_{ij} = e_{ij}, \qquad j = 4, 9, 10.$$

The response $Y_i$ is computed as

$$Y_i = \sum_{j=1}^{10} X_{ij}\beta_{ij} + g(T_i) + \varepsilon_i, \qquad i = 1, \ldots, 100,$$

where $\beta_j = j$, $1 \leq j \leq 4$, $\beta_j = 0$, $5 \leq j \leq 10$, and $\varepsilon_i$'s are sampled from $N(0,1)$. For each $\rho = 0, 0.2, 0.5, 0.8$, we generated $N = 100$ data sets. For comparison, we apply the SCAD penalized regression method, treating $T_i$ as a linear predictor like $X_{ij}$'s. The corresponding estimator is abbreviated as LS-SCAD estimator. Also, profile least squares without variable selection (PLM), profile least squares using AIC for variable selection (PLM-AIC) and partially linear regression using Lasso (PLM-LASSO) are applied for comparison. We investigate two different $g(\cdot)$ functions: Scenario 1, $g(t) = \cos(t)$, and Scenario 2, $g(t) = \cos(2\pi t)$.

The results are summarized in Tables 1 and 2. Columns 3–6 in Table 1 are the averages of the estimates of $\beta_j$, $j = 1, \ldots, 4$, respectively. Column 7 is the number of estimates of $\beta_j$, $5 \le j \le 10$, that are 0, averaged over 100 simulations, and their medians are given in column 8. Column 9 only makes sense for the LS-SCAD estimator. It gives the percentage of times in the 100 simulations in which the coefficient estimate of $T$ equals 0. Model errors are computed as $(\widehat{\boldsymbol{\beta}} - \boldsymbol{\beta})' \operatorname{Cov}(\mathbf{X})(\widehat{\boldsymbol{\beta}} - \boldsymbol{\beta})$. Their medians are listed in the last column, followed by the model errors' standard deviations in parentheses.

In Scenario 1, the nonparametric part $g(T) = \cos(T)$ can be fairly well approximated by a linear function on $[0, 1]$. As a result, the LS-SCAD estimator is expected to give good estimates. It is shown in Table 1 that the estimates of $\beta_j, 1 \le j \le 4$, are all very close to the underlying values. LS-SCAD and PLM-SCAD pick out the covariates with zero coefficients efficiently. LS-SCAD has similar performance to the PLM-AIC and PLM-SCAD in variable selection. On average, each time 83% of the covariates with zero coefficients are selected, and none of the covariates with nonzero coefficients are incorrectly chosen as trivial in the 100 simulations. PLM-LASSO already significantly shrinks the estimates before detecting all the coefficients equal to 0. In each design setting, about 2/3 of the time, the LS-SCAD method attributes no effect to $T$, which does have a quasi-linear effect on $Y$. This is due to the relatively small variation caused in $g(T)$ (with a range less than 0.5), compared with the random variation. Despite this, it performs best with respect to the model error associated with the $\mathbf{X}$ part. PLM-SCAD outperforms PLM-AIC in model errors and is more competent than PLM-LASSO in variable selection.

In Scenario 2, $g(T) = \cos(2\pi T)$. This change in $g(T)$ makes it hard to have a linear approximation of $g(T)$ on $[0, 1]$. So, the LS-SCAD estimator is expected to fail in this situation. Besides, the variation in $g(\cdot)$ (with a range of 2) is relatively large compared to the variation in the error term. Thus, misspecification of $g(T)$ introduces bias in estimating $\boldsymbol{\beta}$. In columns 3–6, with respect to the LS-SCAD estimator, the estimates of the nonzero coefficients are clearly biased, and the biases become larger as the correlation between covariates increases.

Table 2 summarizes the performance of the sandwich estimator of the standard error of the PLM-SCAD estimator for Scenario 1. Columns 2, 4, 6 and 8 are the standard errors of $\beta_j, 1 \le j \le 4$, in the 100 simulations, respectively, while columns 3, 5, 7 and 9 are the average of the standard deviation estimates of these coefficients, obtained via the Hessian matrices. It is seen that the sandwich estimator of the standard error works well, though it slightly underestimates the sampling variation.

TABLE 1
*Example 1, comparison of estimators*

| Estimator | $\rho$ | $\beta_1$ | $\beta_2$ | $\beta_3$ | $\beta_4$ | $\overline{K}$ | $\widetilde{K}$ | %($\widehat{g}(T)=0$) | MME($\times 10^{-2}$) (SD) |
|---|---|---|---|---|---|---|---|---|---|
| Scenario 1 | | | | | | | | | |
| LS-SCAD | 0 | 0.982 | 2.006 | 2.981 | 4.002 | 4.66 | 5 | 70 | 5.48 (4.41) |
| | 0.2 | 1.003 | 1.997 | 2.980 | 3.998 | 4.36 | 5 | 53 | 6.95 (5.11) |
| | 0.5 | 1.012 | 1.997 | 2.974 | 4.014 | 4.40 | 5 | 57 | 7.63 (6.24) |
| | 0.8 | 1.034 | 2.004 | 2.938 | 4.028 | 4.61 | 5 | 66 | 9.74 (10.23) |
| PLM | 0 | 0.988 | 1.984 | 3.001 | 4.005 | 0 | 0 | 0 | 11.62 (5.91) |
| | 0.2 | 1.015 | 1.970 | 3.004 | 3.998 | 0 | 0 | 0 | 12.03 (5.69) |
| | 0.5 | 1.026 | 1.965 | 3.005 | 3.996 | 0 | 0 | 0 | 12.29 (6.47) |
| | 0.8 | 1.048 | 1.948 | 3.008 | 3.994 | 0 | 0 | 0 | 15.19 (11.40) |
| PLM-AIC | 0 | 0.989 | 1.984 | 2.998 | 4.005 | 4.89 | 5 | 0 | 9.22 (5.89) |
| | 0.2 | 1.019 | 1.969 | 3.007 | 3.998 | 4.75 | 5 | 0 | 9.44 (5.88) |
| | 0.5 | 1.027 | 1.964 | 3.009 | 4.001 | 4.78 | 5 | 0 | 10.32 (6.45) |
| | 0.8 | 1.046 | 1.949 | 3.009 | 4.005 | 4.81 | 5 | 0 | 13.19 (10.78) |
| PLM-LASSO | 0 | 0.937 | 1.934 | 2.946 | 3.952 | 2.42 | 2 | 0 | 8.00 (5.89) |
| | 0.2 | 0.971 | 1.935 | 2.971 | 3.951 | 2.48 | 2 | 0 | 8.98 (5.31) |
| | 0.5 | 0.991 | 1.947 | 2.992 | 3.951 | 2.92 | 3 | 0 | 8.44 (5.68) |
| | 0.8 | 1.045 | 1.939 | 3.011 | 3.928 | 3.56 | 4 | 0 | 10.52 (10.58) |
| PLM-SCAD | 0 | 0.988 | 1.985 | 2.999 | 4.004 | 4.49 | 5 | 0 | 6.27 (5.50) |
| | 0.2 | 1.017 | 1.970 | 3.007 | 3.997 | 4.46 | 5 | 0 | 6.78 (5.11) |
| | 0.5 | 1.027 | 1.965 | 3.009 | 3.998 | 4.69 | 5 | 0 | 7.56 (6.07) |
| | 0.8 | 1.045 | 1.948 | 3.014 | 4.001 | 4.78 | 5 | 0 | 12.33 (11.03) |
| Scenario 2 | | | | | | | | | |
| LS-SCAD | 0 | 0.923 | 2.147 | 3.066 | 3.997 | 4.55 | 5 | 67 | 14.53 (9.47) |
| | 0.2 | 0.925 | 2.145 | 3.033 | 3.968 | 4.55 | 5 | 72 | 12.05 (10.72) |
| | 0.5 | 0.857 | 2.216 | 3.005 | 3.935 | 4.43 | 5 | 58 | 15.74 (17.65) |
| | 0.8 | 0.606 | 2.559 | 2.916 | 3.871 | 4.64 | 5 | 23 | 59.70 (59.16) |
| PLM | 0 | 0.988 | 1.984 | 3.001 | 4.005 | 0 | 0 | 0 | 11.65 (5.93) |
| | 0.2 | 1.015 | 1.970 | 3.004 | 3.998 | 0 | 0 | 0 | 11.97 (5.70) |
| | 0.5 | 1.026 | 1.965 | 3.005 | 3.996 | 0 | 0 | 0 | 12.31 (6.49) |
| | 0.8 | 1.048 | 1.948 | 3.008 | 3.994 | 0 | 0 | 0 | 15.17 (11.41) |
| PLM-AIC | 0 | 0.989 | 1.984 | 2.998 | 4.005 | 4.89 | 5 | 0 | 9.39 (5.92) |
| | 0.2 | 1.019 | 1.969 | 3.007 | 3.998 | 4.75 | 5 | 0 | 9.44 (5.93) |
| | 0.5 | 1.027 | 1.964 | 3.009 | 4.001 | 4.78 | 5 | 0 | 10.17 (6.48) |
| | 0.8 | 1.046 | 1.949 | 3.009 | 4.005 | 4.81 | 5 | 0 | 13.17 (10.75) |
| PLM-LASSO | 0 | 0.937 | 1.934 | 2.946 | 3.952 | 2.44 | 2.5 | 0 | 7.99 (5.91) |
| | 0.2 | 0.971 | 1.935 | 2.971 | 3.951 | 2.46 | 3 | 0 | 8.93 (5.32) |
| | 0.5 | 0.991 | 1.947 | 2.992 | 3.951 | 2.90 | 3 | 0 | 8.44 (5.72) |
| | 0.8 | 1.045 | 1.939 | 3.011 | 3.928 | 3.53 | 4 | 0 | 10.62 (10.56) |
| PLM-SCAD | 0 | 0.988 | 1.985 | 2.999 | 4.004 | 4.49 | 5 | 0 | 6.29 (5.54) |
| | 0.2 | 1.017 | 1.970 | 3.007 | 3.997 | 4.46 | 5 | 0 | 6.76 (5.10) |
| | 0.5 | 1.027 | 1.965 | 3.009 | 3.998 | 4.69 | 5 | 0 | 7.57 (6.09) |
| | 0.8 | 1.045 | 1.948 | 3.014 | 4.001 | 4.78 | 5 | 0 | 12.70 (11.02) |

TABLE 2
*Example 1, standard errors of the PLM-SCAD estimates*

| $\rho$ | $SD(\beta_1)$ | $\widehat{se}(\beta_1)$ | $SD(\beta_2)$ | $\widehat{se}(\beta_2)$ | $SD(\beta_3)$ | $\widehat{se}(\beta_3)$ | $SD(\beta_4)$ | $\widehat{se}(\beta_4)$ |
|---|---|---|---|---|---|---|---|---|
| 0   | 0.0954 | 0.0947 | 0.1062 | 0.0975 | 0.0988 | 0.0970 | 0.1058 | 0.0966 |
| 0.2 | 0.1003 | 0.0965 | 0.0911 | 0.1007 | 0.1067 | 0.0993 | 0.1113 | 0.0991 |
| 0.5 | 0.1132 | 0.1102 | 0.1128 | 0.1249 | 0.1264 | 0.1248 | 0.1318 | 0.1139 |
| 0.8 | 0.1590 | 0.1608 | 0.1922 | 0.2062 | 0.2024 | 0.2073 | 0.2116 | 0.1773 |

We have also examined the behavior of the PLM-SCAD estimator of $g(\cdot)$ in Scenario 2. The estimator performs well and is globally close to the true curves (plot not shown). In particular, its performance gets better as $\rho$ decreases.

EXAMPLE 2. The PLM-SCAD estimation is implemented in the analysis of the workers' wage data from Berndt (1991). This data set contains the wage information of 534 workers and their education, living region, gender, race, occupation and marriage status information. Also given are their years of experience. It is not appropriate to assume a linear relationship between years of experience and wage level. However, the main concern is how important the other variables are to wage. In particular, we consider

$$Y_i = g(T_i) + \sum_{j=1}^{14} X_{ij}\beta_j + \varepsilon_i, \qquad i = 1, \ldots, 534,$$

where $Y_i$ is the $i$th worker's wage, $T_i$ is his years of experience, $X_{ij}$ is his $j$th variable and $\varepsilon_i$'s are i.i.d variables with mean 0 and finite variance. There are 14 covariates besides the years of experience. Brief description of the variables, as well as the PLM-SCAD estimates of $\beta_j$'s, can be found in Table 3. As a comparison, the estimates of $\beta_j$ from the unpenalized PLM and Lasso-penalized PLM are given in the third and fourth columns, respectively. PLM-SCAD selects 10 of the 14 covariates, while PLM-LASSO keeps 12.

**6. Discussion.** In this paper, we studied the SCAD-penalized method for variable selection and estimation in the PLM with a divergent number of covariates. B-spline basis functions are used for fitting the nonparametric part. Variable selection and coefficient estimation in the parametric part are achieved simultaneously. The oracle property of the PLM-SCAD estimator

TABLE 3
*Wage data example*

| Variable | Description | $\widehat{\beta}$ (SE) | $\widehat{\beta}_{\text{LASSO}}$ (SE) | $\widehat{\beta}_{\text{SCAD}}$ (SE) |
|---|---|---|---|---|
| $X_1$ | Number of years of education | 0.621 (0.102) | 0.616 (0.086) | 0.645 (0.092) |
| $X_2$ | 1 = southern region, 0 = other | −0.451 (0.417) | −0.313 (0.221) | −0.206 (0.153) |
| $X_3$ | 1 = Female, 0 = Male | −1.956 (0.417) | −1.790 (0.334) | −2.010 (0.388) |
| $X_4$ | 1 = union member, 0 = nonmember | 1.602 (0.508) | 1.343 (0.382) | 1.374 (0.419) |
| $X_5$ | 1 = black, 0 = other | −0.869 (0.574) | −0.516 (0.284) | −0.428 (0.233) |
| $X_6$ | 1 = Hispanic, 0 = other | −0.588 (0.868) | −0.088 (0.060) | 0 (−) |
| $X_7$ | 1 = management, 0 = other | 3.433 (0.796) | 2.909 (0.523) | 3.316 (1.021) |
| $X_8$ | 1 = sales, 0 = other | −0.498 (0.855) | −0.311 (0.192) | −0.057 (0.047) |
| $X_9$ | 1 = clerical, 0 = other | 0.149 (0.683) | 0 (−) | 0 (−) |
| $X_{10}$ | 1 = service, 0 = other | −0.468 (0.680) | −0.494 (0.275) | −0.223 (0.148) |
| $X_{11}$ | 1 = professional, 0 = other | 2.143 (0.731) | 1.781 (0.432) | 2.011 (0.526) |
| $X_{12}$ | 1 = manufacturing, 0 = other | 1.162 (0.595) | 0.799 (0.329) | 0.843 (0.278) |
| $X_{13}$ | 1 = construction, 0 = other | 0.678 (0.962) | 0.075 (0.048) | 0 (−) |
| $X_{14}$ | 1 = married, 0 = other | −0.008 (0.421) | 0 (−) | 0 (−) |

*Notes.* Columns 3–5 are the estimates of $\beta_j$, $j = 1, \ldots, 14$. Their corresponding standard errors are given in parentheses following them.

of the parametric part was established, and consistency of the PLM-SCAD estimator of the nonparametric part was shown.

We have focused on the case where there is one variable in the nonparametric part. Nonetheless, this may be extended to the case of $d$ covariates $T_1, \ldots, T_d$. Specifically, consider the model

$$(5) \qquad Y = \mathbf{X}^{(n)\prime}\boldsymbol{\beta}^{(n)} + g(T_1, \ldots, T_d) + \varepsilon.$$

The PLM-SCAD estimator $(\widehat{\boldsymbol{\beta}}^{(n)}, \widehat{g}_n)$ can be obtained via

$$\min_{(\mathbf{b}^{(n)} \in \mathcal{R}^{p_n}, \phi \in \mathcal{S})} \left\{ \sum_{i=1}^{n} (Y_i - \mathbf{X}_i^{(n)\prime}\mathbf{b}^{(n)} - \phi)^2 + n \sum_{j=1}^{p_n} p_{\lambda_n}(b_j^{(n)}; a) \right\}.$$

Here, $\mathcal{S}$ is the space of all the $d$-variate functions on $[0, 1]^d$ that meet some requirement of smoothness. In particular, we can take $\mathcal{S}$ to be the space of the products of the B-spline basis functions, then project $\mathbb{X}^{(n)}$ and $\mathbf{Y}$ onto this space with this basis and perform the SCAD-penalized regression to $\widetilde{\mathbf{Y}}$ on $\widetilde{\mathbb{X}}^{(n)}$. This has already been discussed in Friedman (1991). However, for large $d$ and moderate sample size, even with very small $M_n$, this model may suffer from the "curse of dimensionality."

A more parsimonious extension is the partially linear additive model (PLAM)

$$(6) \qquad Y = \mu + \mathbf{X}^{(n)\prime}\boldsymbol{\beta}^{(n)} + \sum_{l=1}^{d} g_l(T_l) + \varepsilon,$$

where $E[g_l(T_l)] = 0$ holds for $l = 1, \ldots, d$. To estimate $\boldsymbol{\beta}$ and $g_l$, for each $T_l$, we first determine the partition $\Delta_{nl}$. For simplicity, we assume that the numbers of knots are $M_n$ and the mesh sizes are $O(M_n^{-1})$ for all $l$. Suppose that $X$ and $Y$ are centered. The PLAM–SCAD estimator $(\widehat{\boldsymbol{\beta}}^{(n)}, \widehat{g}_{n1}, \ldots, \widehat{g}_{n1})$ is then defined to be the minimizer of

$$\sum_{i=1}^{n}\left[Y_i - \mathbf{X}_i^{(n)\prime}\mathbf{b}^{(n)} - \sum_{l=1}^{d}\phi_l(T_{il})\right]^2 + n\sum_{j=1}^{p_n} p_{\lambda_n}(b_j^{(n)}; a),$$

subject to $\sum_{i=1}^{n} \phi_l(T_{il}) = 0$ and $\phi_l$ is an element of $\mathcal{S}_m(\Delta_{nl})$.

Under the assumptions similar to those for the PLM-SCAD estimator, $\widehat{\boldsymbol{\beta}}^{(n)}$ can be shown to possess the oracle property. Furthermore, if the joint distribution of $(T_1, \ldots, T_d)$ is absolutely continuous and its density is bounded away from 0 and infinity on $[0,1]^d$, following the proof of Lemma 7 in Stone (1985) and that of Theorem 4 here, we can obtain the same global consistency rate for each additive component, that is,

$$\|\widehat{g}_{nl} - g_l\|_{T_l} = O_P(k_n/\sqrt{n} + \sqrt{M_n/n} + \sqrt{k_n}M_n^{-s_g}), \qquad l = 1, \ldots, d.$$

One way to compute the PLAM-SCAD estimator is the following. First, form the B-spline basis $\{B_{nw}, 1 \le w \le q_n\}$ as follows: the first $M_n + m - 1$ components are the B-spline basis functions corresponding to $T_1$ ignoring the intercept, the second $M_n + m - 1$ components corresponding to $T_2$, and so on. The intercept is the last component. So here, $q_n = dM_n + dm - d + 1$. Now computation can proceed in a similar way to that for the PLM-SCAD estimator.

Our results require that $p_n < n$. While this condition is often satisfied in applications, there are important settings in which it is violated. For example, in studies with microarray data as covariate measurements, the number of genes (covariates) is typically greater than the sample size. Without any further assumptions on the structure of covariate matrix, the regression parameter is in general not identifiable if $p_n > n$. It is an interesting topic of future research to identify conditions under which the PLM-SCAD estimator achieves consistent variable selection and asymptotic normality, even when $p_n > n$.

## APPENDIX

Before embarking on proving the asymptotic results, we give an overview of how the proofs are related to those in Fan and Peng (2004). In their work, the subject under study is a local minimizer of the objective function. We look for conditions when the global minimizer enjoys the desirable properties. In the absence of a nonparametric term, Huang, Horowitz and Ma (2008) solved this problem for the bridge estimator. Identifiability of each component in **X** is a basic requirement in both works. The partial residual approach changes the partially linear model to a linear model by smoothing out the effect of $T$. Identifiability requires that none of the components in **X** or any of their linear combinations vanish after smoothing. With a divergent $p$, uniformness among the components of **X** in smoothness is necessary. Once a certain convergence rate is assured, the proofs for Theorems 2 and 3 are similar to their proofs for consistent variable selection and efficient estimation. The proof of Theorem 4 combines the results in Stone (1985) about the convergence rate of nonparametric regression and the oracle property obtained in this paper.

We now give the proofs of the results stated in Section 3. Write

$$\mathbb{X}^{(n)} = (X_{ij})_{\substack{i=1,\ldots,n \\ j=1,\ldots,p_n}} = (\theta_j^{(n)}(T_i))_{\substack{i=1,\ldots,n \\ j=1,\ldots,p_n}} + (e_{ij}^{(n)})_{\substack{i=1,\ldots,n \\ j=1,\ldots,p_n}}$$

$$\triangleq \boldsymbol{\theta}^{(n)}(\mathbf{T}) + \mathbb{E}_n,$$

and $I - P_{\mathbf{Z}}^{(n)}$ is written as **W** for simplicity.

LEMMA 1.  *Under* (A1), $\|\mathbb{X}^{(n)\prime}\mathbf{W}\mathbb{X}^{(n)}/n - \Xi^{(n)}\| \xrightarrow{P} 0$.

PROOF.  For simplicity, write $\mathbf{A}^{(n)} = \mathbb{X}^{(n)\prime}\mathbf{W}\mathbb{X}^{(n)}/n$ and $\mathbf{C}^{(n)} = \mathbf{A}^{(n)} - \Xi^{(n)}$. Note that $\mathbb{X}_{\cdot j}^{(n)} = \mathbf{e}_{\cdot j}^{(n)} + \theta_{nj}(\mathbf{T})$, where $\mathbf{e}_{\cdot j}^{(n)} = (e_{1j}^{(n)}, \ldots, e_{nj}^{(n)})'$.

$$|C_{jl}^{(n)}| = \left| \left( \frac{\mathbf{e}_{\cdot j}^{(n)\prime}\mathbf{e}_{\cdot l}^{(n)}}{n} - \Xi_{jl}^{(n)} \right) + \frac{\mathbf{e}_{\cdot j}^{(n)\prime} P_{\mathbf{Z}}^{(n)} \mathbf{e}_{\cdot l}^{(n)}}{n} + \frac{\mathbf{e}_{\cdot j}^{(n)\prime} \mathbf{W} \boldsymbol{\theta}_l^{(n)}(\mathbf{T})}{n} \right.$$

$$\left. + \frac{\mathbf{e}_{\cdot l}^{(n)\prime} \mathbf{W} \boldsymbol{\theta}_j^{(n)}(\mathbf{T})}{n} + \frac{\boldsymbol{\theta}_j^{(n)}(\mathbf{T})' \mathbf{W} \boldsymbol{\theta}_l^{(n)}(\mathbf{T})}{n} \right|.$$

By Condition 2, $E[n^{-1}\mathbf{e}_{\cdot j}^{(n)\prime}\mathbf{e}_{\cdot l}^{(n)} - \Xi_{jl}^{(n)}]^2 = n^{-1}\operatorname{Var}(e_j^{(n)} e_l^{(n)}) \leq n^{-1}\sigma_{4e}$, since

$$E[n^{-1}\mathbf{e}_{\cdot j}^{(n)\prime} P_{\mathbf{Z}}^{(n)} \mathbf{e}_{\cdot j}^{(n)}]^2 = n^{-2} E\{E[(\mathbf{e}_{\cdot j}^{(n)\prime} P_{\mathbf{Z}}^{(n)} \mathbf{e}_{\cdot j}^{(n)})^2 | \mathbb{Z}^{(n)}]\}$$

$$= n^{-2} E\left\{ \sum_{i=1}^{n} \sum_{i'=1}^{n} \sum_{\iota=1}^{n} \sum_{\iota'=1}^{n} P_{ii'} P_{\iota\iota'} E[e_{ij}^{(n)} e_{i'j}^{(n)} e_{\iota j}^{(n)} e_{\iota' j}^{(n)} | \mathbb{Z}^{(n)}] \right\},$$

and

$$P_{ii'} P_{\iota\iota'} E[e_{ij}^{(n)} e_{i'j}^{(n)} e_{\iota j}^{(n)} e_{\iota' j}^{(n)} | \mathbb{Z}^{(n)}] = \begin{cases} P_{ii} P_{\iota\iota} \Sigma_{jj}^{(n)}(T_i) \Sigma_{jj}^{(n)}(T_\iota), & i = i' \neq \iota = \iota', \\ P_{ii'}^2 \Sigma_{jj}^{(n)}(T_i) \Sigma_{jj}^{(n)}(T_{i'}), & i = \iota \neq i' = \iota', \\ P_{ii'}^2 \Sigma_{jj}^{(n)}(T_i) \Sigma_{jj}^{(n)}(T_{i'}), & i = \iota' \neq i' = \iota, \\ P_{ii}^2 E[e_{ij}^{(n)^4} | T_i], & i = i' = \iota = \iota', \\ 0, & \text{otherwise}, \end{cases}$$

together with $\Sigma_{jj}^{(n)}(T_i) \leq \sigma_{4e}^{1/2}$ and $P_{\mathbf{Z},ii}^{(n)} \leq 1$, we have

$$E[n^{-1}\mathbf{e}_{\cdot j}^{(n)\prime} P_{\mathbf{Z}}^{(n)} \mathbf{e}_{\cdot j}^{(n)}]^2 \leq n^{-2}\sigma_{4e}\{E[\text{tr}^2(P_{\mathbf{Z}}^{(n)})] + 2E[\text{tr}(P_{\mathbf{Z}}^{(n)2})]\}$$
$$+ n^{-2}\sigma_{4e}E[\text{tr}(P_{\mathbf{Z}}^{(n)})]$$
$$\leq n^{-2}\sigma_{4e}(q_n^2 + 3q_n).$$

By Corollary 6.21 in Schumaker (1981) and the properties of least square regression,

$$E[n^{-1}\boldsymbol{\theta}_j^{(n)}(\mathbf{T})\mathbf{W}\boldsymbol{\theta}_j^{(n)}(\mathbf{T})] \leq C_1 M_\theta(\overline{\Delta}_n)^{2s_\theta},$$

where $C_1$ is a constant determined only by $r_\theta$. By the Cauchy–Schwarz inequality and $C_r$ inequality, we have

$$\|\mathbf{C}^{(n)}\|^2 = O_P(p_n^2/n + p_n^2 M_n^2/n^2 + p_n^2 M_n^{-2s_\theta}).$$

The convergence follows from (A1). □

LEMMA 2. $E[\text{tr}(\mathbb{X}^{(n)\prime}\mathbf{W}\mathbb{X}^{(n)})] = O(np_n).$

PROOF. We have

$$E[\text{tr}(\mathbb{X}^{(n)\prime}\mathbf{W}\mathbb{X}^{(n)})]$$
$$= E[\text{tr}([\mathbb{E}^{(n)} + \boldsymbol{\theta}^{(n)}(\mathbf{T})]'\mathbf{W}[\mathbb{E}^{(n)} + \boldsymbol{\theta}^{(n)}(\mathbf{T})])]$$
$$= E[\text{tr}(\mathbb{E}^{(n)\prime}\mathbf{W}\mathbb{E}^{(n)} + 2\mathbb{E}^{(n)\prime}\mathbf{W}\boldsymbol{\theta}^{(n)}(\mathbf{T}) + \boldsymbol{\theta}^{(n)}(\mathbf{T})'\mathbf{W}\boldsymbol{\theta}^{(n)}(\mathbf{T}))]$$
$$= E\{E[\text{tr}(\mathbb{E}^{(n)\prime}\mathbf{W}\mathbb{E}^{(n)} + 2\mathbb{E}^{(n)\prime}\mathbf{W}\boldsymbol{\theta}^{(n)}(\mathbf{T}) + \boldsymbol{\theta}^{(n)}(\mathbf{T})'\mathbf{W}\boldsymbol{\theta}^{(n)}(\mathbf{T}))|\mathbf{T}]\}$$
$$= E\{E[\text{tr}(\mathbb{E}^{(n)\prime}\mathbf{W}\mathbb{E}^{(n)})|\mathbf{T}]\} + E[\text{tr}(\boldsymbol{\theta}^{(n)}(\mathbf{T})'\mathbf{W}\boldsymbol{\theta}^{(n)}(\mathbf{T}))]$$
$$\leq E\{E[\text{tr}(\mathbb{E}^{(n)\prime}\mathbf{W}\mathbb{E}^{(n)})|\mathbf{T}]\} + C_1 n p_n M_\theta M_n^{-2s_\theta}$$
$$= E\left\{E\left[\sum_{j=1}^{p_n} \mathbf{e}_{\cdot j}^{(n)\prime} \mathbf{W} \mathbf{e}_{\cdot j}^{(n)} \Big| \mathbf{T}\right]\right\} + C_1 n p_n M_\theta M_n^{-2s_\theta}$$
$$= E\left[\sum_{j=1}^{p_n} \text{tr}(\mathbf{W}\Sigma_{jj}^{(n)}(\mathbf{T}))\right] + C_1 n p_n M_\theta M_n^{-2s_\theta}$$
$$\leq n p_n \sigma_{4e}^{1/2} + C_1 n p_n M_\theta M_n^{-2s_\theta}, \qquad \text{tr}(AB) \leq \lambda_{\max}(B)\text{tr}(A).$$

Here, $\Sigma_{jj}^{(n)}(\mathbf{T}) = \text{diag}(\Sigma_{jj}^{(n)}(T_1),\ldots,\Sigma_{jj}^{(n)}(T_n))$. □

PROOF OF THEOREM 1. Let $\boldsymbol{\varepsilon} = (\varepsilon_1,\ldots,\varepsilon_n)'$ and $\mathbf{g}(\mathbf{T}) = (g(T_1),\ldots,g(T_n))'$. Since $\widehat{\boldsymbol{\beta}}^{(n)}$ minimizes $Q_n(\mathbf{b}^{(n)})$, it necessarily holds that $Q_n(\widehat{\boldsymbol{\beta}}^{(n)}) \leq$

$Q_n(\boldsymbol{\beta}^{(n)})$. Rewriting this inequality, we have

$$\|\mathbf{W}\mathbb{X}^{(n)}(\widehat{\boldsymbol{\beta}}^{(n)} - \boldsymbol{\beta}^{(n)})\|^2 - 2(\boldsymbol{\varepsilon} + \mathbf{g}(\mathbf{T}))'\mathbf{W}\mathbb{X}^{(n)}(\widehat{\boldsymbol{\beta}}^{(n)} - \boldsymbol{\beta}^{(n)}) \leq \frac{nk_n}{2}(a+1)\lambda_n^2.$$

Let $\boldsymbol{\delta}_n = n^{-1/2}[\mathbb{X}^{(n)'}\mathbf{W}\mathbb{X}^{(n)}]^{1/2}(\widehat{\boldsymbol{\beta}}^{(n)} - \boldsymbol{\beta}^{(n)})$, and

$$\boldsymbol{\omega}_n = n^{-1/2}[\mathbb{X}^{(n)'}\mathbf{W}\mathbb{X}^{(n)}]^{-1/2}\mathbb{X}^{(n)'}\mathbf{W}(\boldsymbol{\varepsilon} + \mathbf{g}(\mathbf{T})).$$

Then, $\|\boldsymbol{\delta}_n - \boldsymbol{\omega}_n\|^2 \leq \|\boldsymbol{\omega}_n\|^2 + 0.5k_n(a+1)\lambda_n^2$. By the $C_r$ inequality,

$$\|\boldsymbol{\delta}_n\|^2 \leq 2(\|\boldsymbol{\delta}_n - \boldsymbol{\omega}_n\|^2 + \|\boldsymbol{\omega}_n\|^2) \leq 4\|\boldsymbol{\omega}_n\|^2 + k_n(a+1)\lambda_n^2.$$

Examine

$$\|\boldsymbol{\omega}_n\|^2 = n^{-1}(\boldsymbol{\varepsilon} + \mathbf{g}(\mathbf{T}))'\mathbf{W}\mathbb{X}^{(n)}[\mathbb{X}^{(n)'}\mathbf{W}\mathbb{X}^{(n)}]^{-1}\mathbb{X}^{(n)'}\mathbf{W}(\boldsymbol{\varepsilon} + \mathbf{g}(\mathbf{T}))$$
$$\triangleq I_{n1} + I_{n2} + I_{n3},$$

where

$$I_{n1} = n^{-1}\boldsymbol{\varepsilon}'\mathbf{W}\mathbb{X}^{(n)}[\mathbb{X}^{(n)'}\mathbf{W}\mathbb{X}^{(n)}]^{-1}\mathbb{X}^{(n)'}\mathbf{W}\boldsymbol{\varepsilon},$$
$$I_{n2} = 2n^{-1}\boldsymbol{\varepsilon}'\mathbf{W}\mathbb{X}^{(n)}[\mathbb{X}^{(n)'}\mathbf{W}\mathbb{X}^{(n)}]^{-1}\mathbb{X}^{(n)'}\mathbf{W}\mathbf{g}(\mathbf{T}),$$
$$I_{n3} = n^{-1}\mathbf{g}(\mathbf{T})'\mathbf{W}\mathbb{X}^{(n)}[\mathbb{X}^{(n)'}\mathbf{W}\mathbb{X}^{(n)}]^{-1}\mathbb{X}^{(n)'}\mathbf{W}\mathbf{g}(\mathbf{T}).$$

Now, $I_{n1} = E[E(I_{n1}|\mathbb{X}^{(n)}, \mathbf{T})]O_P(1) = p_n n^{-1}O_P(1)$. By the property of projection matrices,

$$I_{n3} \leq n^{-1}\mathbf{g}(\mathbf{T})'\mathbf{W}\mathbf{g}(\mathbf{T}) = M_n^{-2s_g}O(1).$$

Thus, $\|\boldsymbol{\omega}_n\|^2 = O_P(p_n/n + M_n^{-2s_g})$. Furthermore,

$$\|\widehat{\boldsymbol{\beta}}^{(n)} - \boldsymbol{\beta}^{(n)}\|^2 = O_P(p_n/n + M_n^{-2s_g} + k_n\lambda_n^2)$$

follows from Lemma 1 with (A2). Thus, (A3) immediately leads to the consistency. □

LEMMA 3 (Rate of convergence). *Suppose* (A1)–(A4) *hold. Then,*

$$\|\widehat{\boldsymbol{\beta}}^{(n)} - \boldsymbol{\beta}^{(n)}\| = O_P(\sqrt{p_n/n} + \sqrt{p_n}/M_n^{s_g}).$$

PROOF. Let $u_n = \sqrt{p_n/n} + M_n^{-s_g} + \sqrt{k_n}\lambda_n$. When $u_n = o(\min_{1\leq j\leq k_n}|\beta_j^{(n)}|)$, with probability tending to 1, $\min_{1\leq j\leq k_n}|\widehat{\beta}_j^{(n)}| > a\lambda_n$. Given a sequence $\{h_n: h_n > 0\}$ that converges to 0, partition $\mathcal{R}^{p_n} \setminus \{\mathbf{0}_{p_n}\}$ into shells $\{S_{n,l}, l = $

$0, 1, \ldots\}$, where $S_{n,l} = \{\mathbf{b}^{(n)} : 2^{l-1} h_n \leq \|\mathbf{b}^{(n)} - \boldsymbol{\beta}^{(n)}\| < 2^l h_n\}$. Then,

$$P(\|\widehat{\boldsymbol{\beta}}_n^{(n)} - \boldsymbol{\beta}^{(n)}\| \geq 2^L h_n)$$

$$\leq o(1) + \sum_{\substack{l > L \\ 2^l h_n \leq 2^{L_1} u_n}} P(\widehat{\boldsymbol{\beta}}_n^{(n)} \in S_{n,l}, \|\mathbf{C}^{(n)}\| \leq c/2)$$

$$\leq o(1) + \sum_{\substack{l > L \\ 2^l h_n \leq 2^{L_1} u_n}} P\left(\inf_{\mathbf{b}^{(n)} \in S_{n,l}} Q_n(\mathbf{b}^{(n)}) \leq Q_n(\boldsymbol{\beta}^{(n)}), \|\mathbf{C}^{(n)}\| \leq \frac{c_\lambda}{2}\right)$$

$$\leq o(1) + \sum_{l > L} P\left(\sup_{\mathbf{b}^{(n)} \in S_{n,l}} 2(\boldsymbol{\varepsilon} + \mathbf{g}(\mathbf{T}))' \mathbf{W} \mathbb{X}^{(n)} (\mathbf{b}^{(n)} - \boldsymbol{\beta}^{(n)})\right.$$

$$\geq \inf_{\mathbf{b}^{(n)} \in S_{n,l}} (\mathbf{b}^{(n)} - \boldsymbol{\beta}^{(n)})' \mathbb{X}^{(n)\prime} \mathbf{W} \mathbb{X}^{(n)} (\mathbf{b}^{(n)} - \boldsymbol{\beta}^{(n)}),$$

$$\left.\|\mathbf{C}^{(n)}\| \leq \frac{c_\lambda}{2}\right)$$

$$\leq \sum_{l > L} P\left(\sup_{\mathbf{b}^{(n)} \in S_{n,l}} (\boldsymbol{\varepsilon} + \mathbf{g}(\mathbf{T}))' \mathbf{W} \mathbb{X}^{(n)} (\mathbf{b}^{(n)} - \boldsymbol{\beta}^{(n)}) \geq 2^{2l-4} n c_\lambda h_n^2\right) + o(1),$$

since

$$E \sup_{\mathbf{b}^{(n)} \in S_{n,l}} |(\boldsymbol{\varepsilon} + \mathbf{g}(\mathbf{T}))' \mathbf{W} \mathbb{X}^{(n)} (\mathbf{b}^{(n)} - \boldsymbol{\beta}^{(n)})|$$

$$\leq 2^l h_n \sqrt{E[(\boldsymbol{\varepsilon} + \mathbf{g}(\mathbf{T}))' \mathbf{W} \mathbb{X}^{(n)} \mathbb{X}^{(n)\prime} \mathbf{W} (\boldsymbol{\varepsilon} + \mathbf{g}(\mathbf{T}))]}$$

$$\leq 2^{l+1/2} h_n \sqrt{E[\boldsymbol{\varepsilon}' \mathbf{W} \mathbb{X}^{(n)} \mathbb{X}^{(n)\prime} \mathbf{W} \boldsymbol{\varepsilon}] + E[\mathbf{g}(\mathbf{T})' \mathbf{W} \mathbb{X}^{(n)} \mathbb{X}^{(n)\prime} \mathbf{W} \mathbf{g}(\mathbf{T})]}$$

$$\leq 2^{l+1/2} h_n \sqrt{C_3 n p_n + E[\mathbf{g}(\mathbf{T})' \mathbf{W} \mathbf{g}(\mathbf{T}) \operatorname{tr}(\mathbb{X}^{(n)} \mathbb{X}^{(n)\prime} \mathbf{W})]}$$

$$\leq 2^l h_n C_4 (\sqrt{n p_n} + n \sqrt{p_n} M_n^{-s_g}).$$

Continuing the previous arguments, by the Markov inequality,

$$P(\|\widehat{\boldsymbol{\beta}}_n^{(n)} - \boldsymbol{\beta}^{(n)}\| \geq 2^L h_n) \leq o(1) + \sum_{l > L} \frac{C_5(\sqrt{p_n} + \sqrt{n p_n} M_n^{-s_g})}{2^{l-4} h_n \sqrt{n}}.$$

This shows that $\|\widehat{\boldsymbol{\beta}}^{(n)} - \boldsymbol{\beta}^{(n)}\| = O_P(\sqrt{p_n/n} + \sqrt{p_n}/M_n^{s_g})$. □

PROOF OF THEOREM 2. Consider the partial derivatives of $Q_n(\boldsymbol{\beta}^{(n)} + \mathbf{v}^{(n)})$. We assume $\|\mathbf{v}^{(n)}\| = O_P(\sqrt{p_n/n} + \sqrt{p_n} M_n^{-s_g})$. Suppose the support sets of $e_j^{(n)}$ are all contained in a compact set $[-C_e, C_e]$. For $j = k_n + 1, \ldots, p_n$, if $\|\mathbf{v}^{(n)}\| \leq \lambda_n$,

$$\frac{\partial Q_n(\boldsymbol{\beta}^{(n)} + \mathbf{v}^{(n)})}{\partial v_j^{(n)}} = 2\mathbb{X}_{\cdot j}^{(n)\prime} \mathbf{W} \mathbb{X}^{(n)} \mathbf{v}^{(n)} + 2\mathbb{X}_{\cdot j}^{(n)\prime} \mathbf{W}(\boldsymbol{\varepsilon} + \mathbf{g}(\mathbf{T})) + n\lambda_n \operatorname{sgn}(v_j^{(n)})$$

$$\triangleq II_{n1,j} + II_{n2,j} + II_{n3,j}.$$

$$\max_{k_n+1 \leq j \leq p_n} |II_{n1,j}| = 2|\mathbb{X}_{\cdot j}^{(n)\prime} \mathbf{W} \mathbb{X}^{(n)} \mathbf{v}^{(n)}|$$

$$\leq 2\|\mathbf{v}^{(n)}\| \max_{k_n+1 \leq j \leq p_n} \|\mathbb{X}_{\cdot j}^{(n)\prime} \mathbf{W} \mathbb{X}^{(n)}\|$$

$$\leq (\sqrt{p_n/n} + \sqrt{p_n} M_n^{-s_g}) O_P(1)$$

$$\times \max_{k_n+1 \leq j \leq p_n} \|\mathbf{W} \mathbb{X}_{\cdot j}^{(n)}\| \lambda_{\max}^{1/2}(\mathbb{X}^{(n)\prime} \mathbf{W} \mathbb{X}^{(n)})$$

$$= (\sqrt{p_n n} + n\sqrt{p_n} M_n^{-s_g}) O_P(1) \sqrt{\lambda_{\max}(\Xi^{(n)})} + o_P(1)$$

$$= \sqrt{p_n(n + n^2 M_n^{-2s_g}) \lambda_{\max}(\Xi^{(n)})} O_P(1).$$

So, this term is dominated by $\frac{1}{2} II_{n3,j}$, as long as

$$\lim \frac{\sqrt{n} \lambda_n}{\sqrt{p_n \lambda_{\max}(\Xi^{(n)})}} = \infty \quad \text{and} \quad \lim \frac{\lambda_n M_n^{s_g}}{\sqrt{p_n \lambda_{\max}(\Xi^{(n)})}} = \infty,$$

both of which are stated in (A5). To sift out all the trivial components, we need

$$P\left(\max_{k_n+1 \leq j \leq p_n} |II_{n2,j}| > n\lambda_n/2\right) \to 0.$$

This is also implied by (A5), as can be seen from

$$P\left(\max_{k_n+1 \leq j \leq p_n} |II_{n2,j}| > n\lambda_n/2\right)$$

$$\leq \frac{2E[\max_{k_n+1 \leq j \leq p_n} |II_{n2,j}|]}{n\lambda_n}$$

$$\leq \frac{2\sqrt{\sum_{j=k_n+1}^{p_n} E[II_{n2,j}^2]}}{n\lambda_n}$$

$$\leq \frac{2\sqrt{2}\sqrt{\sum_{j=k_n+1}^{p_n} \{E[\boldsymbol{\varepsilon}' \mathbf{W} \mathbb{X}_j^{(n)} \mathbb{X}_j^{(n)\prime} \mathbf{W} \boldsymbol{\varepsilon}] + E[g(\mathbf{T})' \mathbf{W} \mathbb{X}_j^{(n)} \mathbb{X}_j^{(n)\prime} \mathbf{W} g(\mathbf{T})]\}}}{n\lambda_n}$$

$$\leq \frac{C_1 \sqrt{nm_n + nM_n^{-2s_g} nm_n}}{n\lambda_n}.$$

The proof is now complete. □

PROOF OF THEOREM 3. Let $\mathbf{A}_n$ be any $\iota \times k_n$ matrix with full row rank and $\Sigma_n = \mathbf{A}_n\mathbf{A}'_n$. From the variable selection conclusion, with probability tending to 1, we have

$$\widehat{\boldsymbol{\beta}}_1^{(n)} - \boldsymbol{\beta}_1^{(n)} = [\mathbb{X}_1^{(n)\prime}\mathbf{W}\mathbb{X}_1^{(n)}]^{-1}\mathbb{X}_1^{(n)\prime}\mathbf{W}(g(\mathbf{T}) + \boldsymbol{\varepsilon}).$$

We consider the limit distribution of

$$\mathbf{V}_n = n^{-1/2}\Sigma_n^{-1/2}\mathbf{A}_n\Xi_{11}^{(n)^{-1/2}}[\mathbb{X}_1^{(n)\prime}\mathbf{W}\mathbb{X}_1^{(n)}](\widehat{\boldsymbol{\beta}}_1^{(n)} - \boldsymbol{\beta}_1^{(n)})$$

$$= n^{-1/2}\Sigma_n^{-1/2}\mathbf{A}_n\Xi_{11}^{(n)^{-1/2}}\mathbb{X}_1^{(n)\prime}\mathbf{W}(g(\mathbf{T}) + \boldsymbol{\varepsilon})$$

$$\triangleq I_{n1} + I_{n2},$$

where

$$I_{n1} = n^{-1/2}\Sigma_n^{-1/2}\mathbf{A}_n\Xi_{11}^{(n)^{-1/2}}\mathbb{X}_1^{(n)\prime}\mathbf{W}g(\mathbf{T})$$

and

$$I_{n2} = n^{-1/2}\Sigma_n^{-1/2}\mathbf{A}_n\Xi_{11}^{(n)^{-1/2}}\mathbb{X}_1^{(n)\prime}\mathbf{W}\boldsymbol{\varepsilon}.$$

Note that the conclusion of Theorem 3 is equivalent to $\mathbf{V}_n \xrightarrow{d} N(\mathbf{0}_\iota, \sigma^2 I_\iota)$. The first term $I_{n1}$ is a $o_P(1)$ term under (A6) and (A7), as shown in

$$I_{n1} = n^{-1/2}\Sigma_n^{-1/2}\mathbf{A}_n\Xi_{11}^{(n)^{-1/2}}\mathbb{E}_1^{(n)\prime}\mathbf{W}g(\mathbf{T})$$

$$+ n^{-1/2}\Sigma_n^{-1/2}\mathbf{A}_n\Xi_{11}^{(n)^{-1/2}}\boldsymbol{\theta}_1^{(n)\prime}(\mathbf{T})\mathbf{W}g(\mathbf{T}),$$

$$= II_{n1} + II_{n2},$$

where

$$\|II_{n1}\|^2 = E\|II_{n1}\|^2 O_P(1)$$

$$= n^{-1}E[g(\mathbf{T})'\mathbf{W}\mathbb{E}_1^{(n)}\Xi_{11}^{(n)^{-1/2}}\mathbf{A}'_n\Sigma_n^{-1}\mathbf{A}_n\Xi_{11}^{(n)^{-1/2}}\mathbb{E}_1^{(n)\prime}\mathbf{W}g(\mathbf{T})]O_P(1)$$

$$= n^{-1}E\{g(\mathbf{T})'\mathbf{W}E[\mathbb{E}_1^{(n)}\Xi_{11}^{(n)^{-1/2}}\mathbf{A}'_n\Sigma_n^{-1}\mathbf{A}_n\Xi_{11}^{(n)^{-1/2}}\mathbb{E}_1^{(n)\prime}|\mathbf{T}]\mathbf{W}g(\mathbf{T})\}$$

$$\times O_P(1)$$

$$\leq n^{-1}E\{g(\mathbf{T})'\mathbf{W}E[\mathbb{E}_1^{(n)}\Xi_{11}^{(n)^{-1}}\mathbb{E}_1^{(n)\prime}|\mathbf{T}]\mathbf{W}g(\mathbf{T})\}O_P(1)$$

$$= n^{-1}E\{g(\mathbf{T})'\mathbf{W}\operatorname{Diag}(\operatorname{tr}(\Xi_{11}^{(n)^{-1}}\Sigma_{11}^{(n)}(T_1)),\ldots,$$

$$\operatorname{tr}(\Xi_{11}^{(n)^{-1}}\Sigma_{11}^{(n)}(T_n)))\mathbf{W}g(\mathbf{T})\}O_P(1)$$

$$\leq n^{-1}\|\mathbf{W}g(\mathbf{T})\|^2 \operatorname{tr}(\Sigma_{u,11}^{(n)})$$

$$= \operatorname{tr}(\Sigma_{u,11}^{(n)})M_n^{-2s_g}O_P(1) = o_P(1)$$

and

$$\|II_{n2}\|^2 \leq n^{-1}\|\mathbf{W}g(\mathbf{T})\|^2 \lambda_{\max}(\mathbf{W}\boldsymbol{\theta}_1^{(n)}(\mathbf{T})\Xi_{11}^{(n)^{-1}}\mathbf{W}\boldsymbol{\theta}_1^{(n)\prime})$$
$$\leq n^{-1}\|\mathbf{W}g(\mathbf{T})\|^2 \|\mathbf{W}\boldsymbol{\theta}_1^{(n)}\|^2 = n^{-1}nM_n^{-2s_g}nM_n^{-2s_\theta}O(1)$$
$$= nM_n^{-2(s_g+s_\theta)}O(1).$$

Decompose the second term $I_{n2}$ as

$$I_{n2} = n^{-1/2}\Sigma_n^{-1/2}\mathbf{A}_n\Xi_{11}^{(n)^{-1/2}}\mathbb{E}_1^{(n)\prime}\boldsymbol{\varepsilon} - n^{-1/2}\Sigma_n^{-1/2}\mathbf{A}_n\Xi_{11}^{(n)^{-1/2}}\mathbb{E}_1^{(n)\prime}P_{\mathbf{Z}}^{(n)}\boldsymbol{\varepsilon}$$
$$+ n^{-1/2}\Sigma_n^{-1/2}\mathbf{A}_n\Xi_{11}^{(n)^{-1/2}}\boldsymbol{\theta}_1^{(n)\prime}(\mathbf{T})\mathbf{W}\boldsymbol{\varepsilon},$$
$$= III_{n1} + III_{n2} + III_{n3}.$$

Actually, the last two terms above are trivial:

$$\|III_{n2}\|^2 = n^{-1}O_P(1)E[\operatorname{tr}(P_{\mathbf{Z}}^{(n)}\mathbb{E}_1^{(n)}\Xi_{11}^{(n)^{-1/2}}\mathbf{A}_n'\Sigma_n^{-1}\mathbf{A}_n\Xi_{11}^{(n)^{-1/2}}\mathbb{E}_1^{(n)\prime}P_{\mathbf{Z}}^{(n)})]$$
$$\leq n^{-1}O_P(1)E[\operatorname{tr}(P_{\mathbf{Z}}^{(n)}\mathbb{E}_1^{(n)}\Xi_{11}^{(n)^{-1}}\mathbb{E}_1^{(n)\prime}P_{\mathbf{Z}}^{(n)})]$$
$$= n^{-1}O_P(1)E[\operatorname{tr}(P_{\mathbf{Z}}^{(n)}\mathbb{E}_1^{(n)}\mathbb{E}_1^{(n)\prime})]$$
$$\leq n^{-1}O_P(1)\operatorname{tr}(\Sigma_{u,11}^{(n)})E[\operatorname{tr}(P_{\mathbf{Z}}^{(n)})]$$
$$= \operatorname{tr}(\Sigma_{u,11}^{(n)})M_n/nO_P(1) = o_P(1).$$

$$\|III_{n3}\|^2 = n^{-1}O_P(1)E[\operatorname{tr}(\mathbf{W}\boldsymbol{\theta}_1^{(n)}(\mathbf{T})\Xi_{11}^{(n)^{-1}}\boldsymbol{\theta}_1^{(n)\prime}(\mathbf{T})\mathbf{W})]$$
$$= k_nM_n^{-2s_\theta}O_P(1) = o_P(1).$$

So we focus on $III_{n1} = n^{-1/2}\Sigma_n^{-1/2}\mathbf{A}_n\Xi_{11}^{(n)^{-1/2}}\mathbb{E}_1^{(n)\prime}\boldsymbol{\varepsilon}$, since

$$\operatorname{Var}(III_{n1}) = E[\operatorname{Var}(III_{n1}|\mathbb{X}^{(n)},\mathbf{T})] = \sigma^2 I_\iota$$

and the infinitely small condition holds, provided $E[\varepsilon^4] < \infty$, and by the Lindeberg–Feller central limit theorem we have $III_{n1} \xrightarrow{d} N(\mathbf{0}_\iota, \sigma^2 I_\iota)$. The conclusion follows from the Slutsky's theorem. □

LEMMA 4. *Sequences of random variables $A_n$ and random vectors $\mathbf{B}_n$ satisfy $E[A_n^2|\mathbf{B}_n] = O_P(u_n^2)$, where $\{u_n\}$ is a sequence of positive numbers. Then, $A_n = O_P(u_n)$.*

PROOF. For any $\varepsilon > 0$, there is some $M_1$, such that $P(E[A_n^2|\mathbf{B}_n] > M_1 u_n^2) < \varepsilon/2$. Let $M_2^2 = 2M_1/\varepsilon$. Then,

$$P(|A_n| > M_2 u_n) \leq P(|A_n| > M_2 u_n, E[A_n^2|\mathbf{B}_n] \leq M_1 u_n^2)$$

$$+ P(E[A_n^2|\mathbf{B}_n] > M_1 u_n^2)$$
$$< E[1_{(|A_n|>M_2 u_n)} 1_{(E[A_n^2|\mathbf{B}_n]\leq M_1 u_n^2)}] + \varepsilon/2$$
$$= E\{1_{(E[A_n^2|\mathbf{B}_n]\leq M_1 u_n^2)} E[1_{(|A_n|>M_2 u_n)}|\mathbf{B}_n]\} + \varepsilon/2$$
$$\leq E\left[1_{(E[A_n^2|\mathbf{B}_n]\leq M_1 u_n^2)} \frac{E[A_n^2|\mathbf{B}_n]}{M_2^2 u_n^2}\right] + \varepsilon/2$$
$$\leq \varepsilon.$$

The arbitrariness of $\varepsilon$ implies the conclusion. $\square$

PROOF OF THEOREM 4. The nonparametric component $g(\cdot)$ at a point $t \in [0,1]$ is estimated with

$$\widehat{g}_n(t) = \mathbf{Z}(t;\Delta_n)'(\mathbb{Z}^{(n)\prime}\mathbb{Z}^{(n)})^{-1}\mathbb{Z}^{(n)\prime}(\mathbf{Y} - \mathbb{X}^{(n)}\widehat{\boldsymbol{\beta}}^{(n)}).$$

With probability tending to 1,

$$\widehat{g}_n(t) - g(t) = \mathbf{Z}(t;\Delta_n)'(\mathbb{Z}^{(n)\prime}\mathbb{Z}^{(n)})^{-1}\mathbb{Z}^{(n)\prime}(\mathbf{Y} - \mathbb{X}_1^{(n)}\widehat{\boldsymbol{\beta}}_1^{(n)}) - g(t)$$
$$= \mathbf{Z}(t;\Delta_n)'(\mathbb{Z}^{(n)\prime}\mathbb{Z}^{(n)})^{-1}\mathbb{Z}^{(n)\prime}g(\mathbf{T}) - g(t)$$
$$+ \mathbf{Z}(t;\Delta_n)'(\mathbb{Z}^{(n)\prime}\mathbb{Z}^{(n)})^{-1}\mathbb{Z}^{(n)\prime}\boldsymbol{\varepsilon}$$
$$- \mathbf{Z}(t;\Delta_n)'(\mathbb{Z}^{(n)\prime}\mathbb{Z}^{(n)})^{-1}\mathbb{Z}^{(n)\prime}\boldsymbol{\theta}_1^{(n)}(\mathbf{T})(\widehat{\boldsymbol{\beta}}_1^{(n)} - \boldsymbol{\beta}_1^{(n)})$$
$$- \mathbf{Z}(t;\Delta_n)'(\mathbb{Z}^{(n)\prime}\mathbb{Z}^{(n)})^{-1}\mathbb{Z}^{(n)\prime}\mathbb{E}_1^{(n)}(\widehat{\boldsymbol{\beta}}_1^{(n)} - \boldsymbol{\beta}_1^{(n)})$$
$$\triangleq I_{n1} + I_{n2} + I_{n3} + I_{n4}.$$

Consider $\|\widehat{g}_n - g\|_T^2 = \int [\widehat{g}_n(t) - g(t)]^2 f_T(t)\,dt$. Without further assumptions, by Lemma 9 in Stone (1985), $\|I_{n1}\|_T^2 = O_P(M_n^{-2s_g})$. When $M_n = o(\sqrt{n})$, by Lemma 4 in Stone (1985),

$$E[\|I_{n2}\|_T^2|\mathbf{T}] = O_P(M_n/n) \quad \text{and hence} \quad \|I_{n2}\|_T^2 = O_P(M_n/n).$$

When $\{\theta_j^{(n)}(\cdot), n \geq 1, 1 \leq j \leq k_n\}$ are uniformly bounded on $[0,1]$,

$$\|I_{n3}\|_T^2 \leq \|\mathbf{Z}(t;\Delta_n)'(\mathbb{Z}^{(n)\prime}\mathbb{Z}^{(n)})^{-1}\mathbb{Z}^{(n)\prime}\boldsymbol{\theta}_1^{(n)}(\mathbf{T})\|_T^2 \|\widehat{\boldsymbol{\beta}}_1^{(n)} - \boldsymbol{\beta}_1^{(n)}\|^2$$
$$\leq [O(k_n) + O_P(k_n M_n^{-2s_\theta})] [O_P(1) M_n^{-2s_g} + k_n/n O_P(1)]$$
$$= O_P(1)(k_n M_n^{-2s_g} + k_n^2 n^{-1}).$$

Similarly,

$$\|I_{n4}\|_T^2 \leq \|\mathbf{Z}(t;\Delta_n)'(\mathbb{Z}^{(n)\prime}\mathbb{Z}^{(n)})^{-1}\mathbb{Z}^{(n)\prime}\mathbb{E}_1^{(n)}\|_T^2 \|\widehat{\boldsymbol{\beta}}_1^{(n)} - \boldsymbol{\beta}_1^{(n)}\|^2$$
$$= \|\widehat{\boldsymbol{\beta}}_1^{(n)} - \boldsymbol{\beta}_1^{(n)}\|^2 \|\mathbf{Z}(t;\Delta_n)'(\mathbb{Z}^{(n)\prime}\mathbb{Z}^{(n)})^{-1}\mathbb{Z}^{(n)\prime}\mathbb{E}_1^{(n)}\|_T^2$$

$$\leq O_P(k_n M_n/n)[O_P(1)M_n^{-2s_g} + k_n/n O_P(1)]$$
$$= O_P(1)(M_n^{1-2s_g} k_n/n + M_n k_n^2/n^2).$$

To sum up, when $k_n = o(\sqrt{n})$, we have $\|\widehat{g}_n - g\|_T^2 = O_P(k_n^2/n + M_n/n + k_n M_n^{-2s_g})$. $\square$

Department of Management Science
University of Miami
Coral Gables, Florida 33124
USA
E-mail: h.xie@miami.edu

Department of Statistics
and Actuarial Science
University of Iowa
Iowa City, Iowa 52242
USA
E-mail: jian-huang@uiowa.edu